\documentclass[12pt]{article}
\usepackage{graphicx} %use graph format
\usepackage{epstopdf}
\usepackage{subfigure}
\usepackage{color}
\usepackage[english]{babel}
\usepackage{amsmath,amsthm}
\usepackage{amssymb}%Áí¼ÓµÄ£¬ÎªÁË´ò³öÕÅÁ¿»ýµÄ·ûºÅ
\usepackage{amsfonts}
\usepackage{booktabs}
\usepackage[colorlinks,linkcolor=red,anchorcolor=blue,citecolor=black]{hyperref} % ¿ÉÒÔÌí¼ÓÒýÓÃ¸ñÊ½

\usepackage{graphicx}
\usepackage{subfigure}
\usepackage{tabularx}
\newtheorem{theorem}{Theorem}
\newtheorem{corollary}{Corollary}

\newtheorem{thm}{Theorem}[section]
\newtheorem{defi}{Definition}[section]
\newtheorem{lem}[thm]{Lemma}

\theoremstyle{plain}

\title{The chromatic-AT choosability of Pseudo-Halin graphs}
\author { Zhiguo Li{$^*$}, Cuicui Ma, Zeling Shao\\
{\small School of Science, Hebei University of Technology, Tianjin 300401, China}\\
{\small Email: zhiguolee@hebut.edu.cn; Zelingshao@163.com; m15102673812@163.com}\\
\date{}
}

\begin{document}
\baselineskip 0.65cm

\maketitle

\begin{abstract}
  The \emph{Alon-Tarsi number} $AT(G)$ of a graph $G$ is the smallest integer $k$ such that $G$ admits an orientation $D$ with maximum outdegree at most $k-1$ in which the number of even Eulerian subdigraphs differs from the number of odd Eulerian subdigraphs. A 2-connected planar graph $G$ with minimum degree at least 3 is called a \emph{pseudo-Halin graph} if deleting the edges on the boundary of a single face $f_0$ yields a tree. Let $\mathcal{A}$ denote the class of pseudo-Halin graphs in which no two interior vertices are adjacent and every interior vertex has odd degree (equivalently, graphs obtained by iterated Haj\'{o}s sums of wheels with an odd number of exterior vertices). In this paper, we completely determine the chromatic number and the Alon-Tarsi number of all pseudo-Halin graphs. Specifically, by structural induction and the method of $AT$-orientations, we prove that $AT(G)=\chi(G)=4$ if $G\in\mathcal{A}$, and $AT(G)=\chi(G)=3$ otherwise. Consequently, every pseudo-Halin graph is chromatic-AT choosable.


\end{abstract}

\section{Introduction}

All graphs considered in this paper are finite, simple, and undirected. Let $G=(V, E)$ be a graph with vertex set $V$ and edge set $E$. The \emph{chromatic number} $\chi(G)$ of $G$ is the minimum number of colors $k$ required to assign colors to the vertices of $G$ such that adjacent vertices receive different colors. By the celebrated Four Color Theorem, $\chi(G)\leq 4$ for every planar graph $G$. A \emph{list assignment} of $G$ is a function $L$ that assigns to each vertex $v\in V$ a set $L(v)$ of permissible colors. A \emph{proper $L$-coloring} of $G$ is a proper vertex coloring $f$ such that $f(v)\in L(v)$ for every $v\in V$. The \emph{choice number} (or \emph{list chromatic number}) $ch(G)$ is the least integer $k$ such that $G$ admits a proper $L$-coloring whenever $|L(v)|\geq k$ for every $v\in V$. Clearly, $\chi(G)\leq ch(G)$ holds for all graphs $G$, and $G$ is said to be \emph{chromatic-choosable} if $\chi(G)=ch(G)$.

In 1992, Alon and Tarsi $\cite{NM1992}$ introduced an algebraic method that relates list colorings of a graph $G$ to Eulerian subdigraphs of its orientations. The \emph{Alon-Tarsi number} of $G$, denoted by $AT(G)$, is the smallest integer $k$ such that there exists an orientation $D$ of $G$ with maximum outdegree at most $k-1$ for which the number of Eulerian subdigraphs with an even number of arcs is different from the number of Eulerian subdigraphs with an odd number of arcs. Alon and Tarsi $\cite{NM1992,TR1995}$ proved that $\chi(G)\leq ch(G)\leq AT(G)$ for every graph $G$. Analogous to chromatic-choosable graphs, a graph $G$ is called \emph{chromatic-AT choosable} if $\chi(G)=AT(G)$.

The Alon-Tarsi number has been extensively investigated for planar graphs and related graph classes. In $\cite{JD2010}$ and $\cite{ESG2023}$, it was shown that for any graph $G$, $\frac{{\rm mad}(G)}{2}< AT(G)\leq \Delta(G) +1$, where ${\rm mad}(G)=\max\{\frac{2|E(H)|}{|V(H)|}: H\subseteq G\}$ is the maximum average degree of $G$. Extending Thomassen's famous 5-choosability theorem, X. Zhu $\cite{XD2017}$ proved that $AT(G)\leq 5$ for every planar graph $G$. Furthermore, Grytczuk and X. Zhu $\cite{JXD2020}$ showed that every planar graph $G$ contains a matching $M$ such that $AT(G-M)\leq 4$, and Kim, Kim, and X. Zhu $\cite{RS2019}$ proved that every planar graph $G$ contains a forest $F$ such that $AT(G-E(F))\leq 3$. For non-planar graphs, T. Abe et al. $\cite{TAK2022}$ proved that every $K_5$-minor-free graph $G$ satisfies $AT(G)\leq 5$. Exact values of the Alon-Tarsi number have also been determined for several structured graph families: Z. Li, Z. Shao, F. Petrov, and A. Gordeev $\cite{ZG2023}$ proved via the polynomial method that for the Cartesian product of two cycles, $AT(C_m\Box C_n)=4$ when $m,n$ are both odd and $AT(C_m\Box C_n)=3$ otherwise; Z. Li, Q. Ye, and Z. Shao $\cite{YQ1,YQ2}$ determined the Alon-Tarsi number of a class of 4-regular planar graphs $R_n$, biwheels $B_n$, and the Mycielski graphs of even cycles, trees, and wheels. In particular, for Halin graphs, Z. Li, Q. Ye, and Z. Shao $\cite{ZG2021}$ proved that $AT(H)=4$ when $H$ is a wheel with an odd number of exterior vertices, and $AT(H)=3$ otherwise.

As a natural generalization of Halin graphs, Liu and Zhang $\cite{LZ2001}$ and Meng, Guo, and Su $\cite{XY2012}$ studied the class of \emph{pseudo-Halin graphs}. A 2-connected planar graph $G$ with minimum degree at least 3 is called a \emph{pseudo-Halin graph} if deleting the edges on the boundary of a single face $f_0$ yields a spanning tree $T$. The face $f_0$ is called the \emph{exterior face} and the other faces are \emph{interior faces}; the vertices on the boundary of $f_0$ are \emph{exterior vertices} and the remaining vertices are \emph{interior vertices}. An exterior vertex $v\in V(f_0)$ is called a \emph{regular vertex} if $d_G(v)=3$, and an \emph{irregular vertex} if $d_G(v)\geq 4$. We denote the sets of regular and irregular exterior vertices by $R(f_0)$ and $IR(f_0)$, respectively. Note that a pseudo-Halin graph is a Halin graph if and only if $IR(f_0)=\emptyset$. When $IR(f_0)\neq \emptyset$, exterior vertices in $IR(f_0)$ have degree at least 2 in the tree $T=G-E(f_0)$, which gives rise to a broader class of graphs than classical Halin graphs.

Equivalently, a pseudo-Halin graph can be written as $G=T\cup C_n$, where $T$ is a spanning tree and $C_n$ is the boundary cycle of the exterior face $f_0$. Let $E^{\prime}=E(G-V(C_n))=\{uv\in E(T)\mid u,v\in V(G)\setminus V(C_n)\}$ be the set of edges in $T$ whose endpoints are both interior vertices. Depending on whether $E^{\prime}=\emptyset$ or $E^{\prime}\neq \emptyset$, pseudo-Halin graphs exhibit different structural properties: when $E^{\prime}=\emptyset$, no two interior vertices are adjacent, and $G$ can be decomposed into an iterated Haj\'{o}s sum of wheels. Let $\mathcal{A}$ be the family of pseudo-Halin graphs $G$ such that $E^{\prime}=\emptyset$ and all interior vertices of $G$ have odd degrees.

Motivated by the results on Halin graphs in $\cite{ZG2021}$ and the chromatic properties of Halin graphs in $\cite{HX1994}$, this paper investigates both the chromatic number and the Alon-Tarsi number of pseudo-Halin graphs. First, we completely determine the chromatic number of all pseudo-Halin graphs:

\begin{theorem}
If $G=T\cup C_n$ is a pseudo-Halin graph, then
$$\chi(G)=\left\{\begin{array}{l}
4,\ \ \ \mbox{\rm if}\ G\in \mathcal{A};\\
3,\ \ \ \mbox{\rm otherwise}.\\
\end{array} \right.$$
\end{theorem}

Furthermore, by analyzing the degeneracy of $G$ for $G\in \mathcal{A}$ and constructing explicit Alon-Tarsi orientations via structural induction for $G\notin \mathcal{A}$, we completely determine the Alon-Tarsi number of all pseudo-Halin graphs in the following unified theorem:

\begin{theorem}
If $G=T\cup C_n$ is a pseudo-Halin graph, then
$$AT(G)=\left\{\begin{array}{l}
4,\ \ \ \mbox{\rm if}\ G\in \mathcal{A};\\
3,\ \ \ \mbox{\rm otherwise}.\\
\end{array} \right.$$
\end{theorem}

Combining Theorems 1 and 2 with the inequality $\chi(G)\leq ch(G)\leq AT(G)$, we immediately obtain the following corollary, which generalizes the corresponding result for Halin graphs in $\cite{ZG2021}$:

\begin{corollary}
Every pseudo-Halin graph $G$ satisfies $\chi(G)=ch(G)=AT(G)$, and hence is chromatic-AT choosable.
\end{corollary}

The remainder of this paper is organized as follows. In Section 2, we present basic definitions and structural lemmas for pseudo-Halin graphs. In Section 3, we prove Theorem 1 by determining the chromatic number of pseudo-Halin graphs in the cases $E^{\prime}\neq \emptyset$ and $E^{\prime}=\emptyset$. In Section 4, we divide the proof of Theorem 2 into two parts corresponding to $G\in \mathcal{A}$ and $G\notin \mathcal{A}$. 

\section{Preliminaries}

\begin{defi} {$\cite{NM1992}$}
A subdigraph $H$ of a directed graph $D$ is called \emph{Eulerian} if $V(H)=V(D)$ and the indegree $d^{-}_{H}(v)$ of every vertex $v$ in $H$ is equal to its outdegree $d^{+}_{H}(v)$. Note that $H$ need not be connected. For a digraph $D$, we denote by $\mathcal{E}(D)$ the family of Eulerian subdigraphs of $D$. An Eulerian subdigraph $H$ is \emph{even} if it has an even number of arcs; otherwise, it is \emph{odd}. Let $\mathcal{E}_{e}(D)$ and $\mathcal{E}_{o}(D)$ denote the families of even and odd Eulerian subdigraphs of $D$, respectively, and let ${\rm diff}(D)=|\mathcal{E}_{e}(D)|-|\mathcal{E}_{o}(D)|$. We say that $D$ is \emph{Alon-Tarsi} if ${\rm diff}(D)\neq0$. If an orientation $D$ of $G$ yields an Alon-Tarsi digraph, then $D$ is called an \emph{Alon-Tarsi orientation} $($or an \emph{$AT$-orientation}, for short$)$ of $G$.
\end{defi}

In general, any orientation $D$ contains at least one even Eulerian subdigraph, namely the empty subdigraph. If an orientation $D$ has no directed cycle, no odd directed cycle, or satisfies $|\mathcal{E}(D)|\equiv 1~({\rm mod}~2)$, it follows readily that ${\rm diff}(D)\neq0$ and hence $D$ is an $AT$-orientation.

\begin{defi} {$\cite{ESG2023}$}
The \emph{degeneracy} of a graph $G$ is the least integer $k$ such that every subgraph $H$ of $G$ has a vertex whose degree in $H$ is at most $k$, and is given by $\max\{\delta(H): H \text{ is a subgraph of } G\}$, where $\delta(H)$ is the minimum degree of the subgraph $H$.
\end{defi}

\begin{lem}{$\cite{ZG2021}$}
If the degeneracy of a graph $G$ is $k$, then $AT(G)\leq k+1$.
\end{lem}

The upper bound $k+1$ in Lemma 2.1 is sharp. For example, the degeneracy of a complete graph $K_n$ is $n-1$, and $K_n$ is chromatic-AT choosable with $AT(K_n)=n$.

\begin{lem}{$\cite{LZ2001}$}
Let $G$ be a pseudo-Halin graph that is not a wheel, with exterior face $f_0$, and let $P=u_1u_2\cdots u_k$ be a longest path in $G-E(f_0)$. Take $w\in \{u_2, u_{k-1}\}$ and let $N(w)$ be the set of neighbours of $w$ in $G$. Then one of the following holds:

$(1)$ The vertex $w$ is an interior vertex of $G$ with $N(w)\subset V(f_0)$ and $|N(w)\cap IR(f_0)|=1$. Let $N(w)=\{y, v_1, v_2, \cdots, v_m\}$ $(m\geq 2)$, where $xv_1$, $yv_m$, $v_iv_{i+1}\in E(f_0)$ $(i=1, 2, \cdots, m-1)$, $y\in IR(f_0)$, $x\neq v_2$, and $y\neq v_{m-1}$. Then the graphs
$$G_1^1:=G-\{w, v_i \mid i=1, 2, \cdots ,m\}+\{xy\}$$
and
$$G_1^2:=G-\{v_i, v_{i+1}, \cdots ,v_j\}+\{v_{i-1}v_{j+1}\}\quad (2\leq i\leq j\leq m-1)$$
are also pseudo-Halin graphs.

$(2)$ The vertex $w$ is an interior vertex of $G$ with $|N(w)\cap (V(G)\setminus V(f_0))|=1$ and $|N(w)\cap R(f_0)|=d(w)-1$. Let $N(w)=\{u, v_1, v_2, \cdots, v_m\}$, where $u$ is an interior vertex, $v_i\in R(f_0)$ $(i=1, 2, \cdots ,m)$, $xv_1$, $yv_m$, $v_iv_{i+1}\in E(f_0)$ $(i=1, 2, \cdots ,m-1)$, $x\neq v_2$, and $y\neq v_{m-1}$. Then the graphs
$$G_2^1:=G-\{v_1, v_2, \cdots, v_m\}+\{xw, yw\}$$
and
$$G_2^2:=G-\{v_i, v_{i+1}, \cdots ,v_j\}+\{v_{i-1}v_{j+1}\}\quad (2\leq i<j\leq m,\ m\geq3)$$
are also pseudo-Halin graphs.
\end{lem}

\begin{defi} {$\cite{NS2023}$}
Given two graphs $G_1=(V_1, E_1)$ and $G_2=(V_2, E_2)$ with disjoint vertex sets, an edge $x_1y_1\in E_1$, and an edge $x_2y_2\in E_2$, the \emph{Haj\'{o}s sum} $H(G_1(x_1y_1), G_2(x_2y_2))=G_1(x_1y_1)+_H G_2(x_2y_2)$ is the graph obtained from the disjoint union $(V_1\cup V_2, E_1\cup E_2)$ by deleting the edges $x_1y_1$ and $x_2y_2$, identifying the vertices $x_1$ and $x_2$ as $v_H(x_1x_2)$, and adding the edge $y_1y_2$.

\begin{figure}[htbp]
\centering
\includegraphics[height=3.5cm, width=0.6\textwidth]{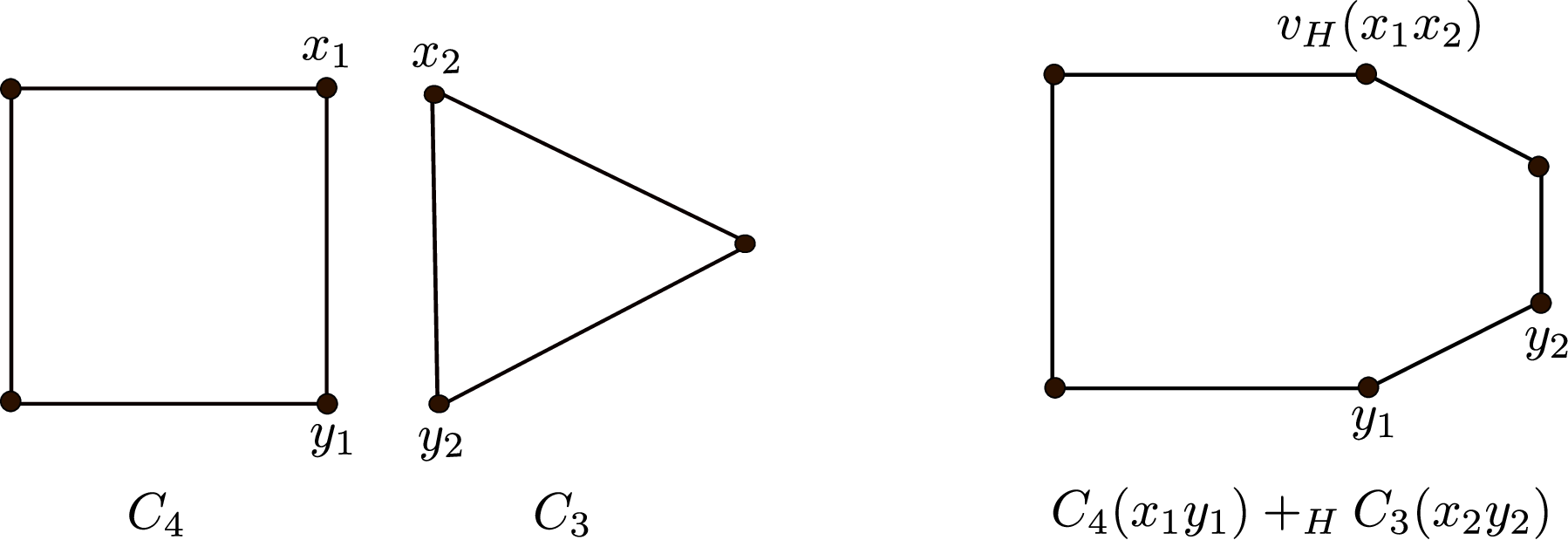}
\centerline{Fig.1  ~The Haj\'{o}s sum of $C_4$ and $C_3$.}
\end{figure}
\end{defi}

\begin{lem}
For every pseudo-Halin graph $G=T\cup C_n$ with $E^{\prime}= \emptyset$, $$G=W_{p_{r+1}}+_HW_{p_{r}}+_H\cdots +_HW_{p_2}+_HW_{p_1},$$ where $W_{p_i}$ is a wheel, $i=1, 2, \cdots, r+1$ $(r\geq 0)$.
\end{lem}

\begin{proof}
The wheel $W_{p_i}$ can be written as $S_{p_i+1}\cup C_{p_i}$, where $S_{p_i+1}$ is the star with $p_i+1$ vertices and $C_{p_i}$ is the cycle with $p_i$ vertices. If $G$ is a wheel, the result holds trivially with $r=0$. Otherwise, for any pseudo-Halin graph $G$ with $E^\prime=\emptyset$, let $P=u_1u_2\cdots u_k$ be a longest path of $T$, and take $w\in \{u_2, u_{k-1}\}$ with $N(w)=\{y, v_1, v_2, \cdots, v_m\}$ $(m\geq2)$ being the set of neighbours of $w$. Since $E^{\prime}=\emptyset$, there is a unique irregular vertex in $N(w)$ and $N(w)\subset V(C_n)$. Assume that $y$ is the irregular vertex and $xv_1$, $yv_m$, $v_iv_{i+1}\in E(C_n)$, $x\neq v_2$, $y\neq v_{m-1}$, $i\in \{1, 2, \cdots, m-1\}$.
Let $R_1=G$, and let $R_2$ be obtained from $R_1$ by identifying $w, v_1, v_2, \cdots, v_m$ with $y$ (equivalently, deleting $w, v_1, \cdots, v_m$ and adding the edge $xy$). By Lemma 2.2, $R_2$ is a pseudo-Halin graph with $E^{\prime}= \emptyset$, and $R_1=R_{2}(yx)+_HW_{p_1}(x_2y_2)$, where $x_2y_2\in E(C_{p_1})$ and $p_1=m+1$. Likewise, $R_2=R_3+_HW_{p_2}$, $\cdots$, $R_{r-1}=R_r+_HW_{p_{r-1}}$. Repeat this process until $R_{r}=W_{p_{r+1}}+_HW_{p_{r}}$. In this process, each graph $R_i$ is a pseudo-Halin graph with $E^{\prime}= \emptyset$, $i=1, 2, \cdots, r$. Therefore $G=W_{p_{r+1}}+_HW_{p_{r}}+_H\cdots +_HW_{p_2}+_HW_{p_1}$.
\end{proof}

\section{Proof of Theorem 1}

It is known in $\cite{HX1994}$ that the chromatic number of a Halin graph is 4 when it is a wheel with an odd number of exterior vertices and 3 otherwise. In this section, we determine the chromatic number of a general pseudo-Halin graph $G=T\cup C_n$ by considering the cases $E^{\prime}\neq \emptyset$ and $E^{\prime}=\emptyset$ in Lemmas 3.1 and 3.2, respectively, thereby proving Theorem 1. Let the exterior vertices $v_1$, $v_2$, $\cdots$, $v_n$ on $C_n$ be listed in counterclockwise order.

\begin{lem}
If $G=T\cup C_n$ is a pseudo-Halin graph with $E^{\prime}\neq \emptyset$, then $\chi(G)=3$.
\end{lem}

\begin{proof}
Denote by $U$ the set of interior vertices that are adjacent to vertices on $C_n$. Obviously, the graph $G$ contains triangles as its subgraphs, so $\chi(G)\geq3$. What is left is to show that $\chi(G)\leq3$.
Since $T$ is a tree, $\chi(T)=2$. Let $\sigma$ be a proper 2-coloring of $T$ using two colors $\{0, 1\}$. Next, a proper 3-coloring of $G$ is constructed by discussing two cases.

{\textbf{Case 1.}}\quad $n$ is even.

Based on $\sigma$, we construct the coloring $\sigma^*$ of $G$ by changing the colors of $v_2$, $v_4$, $\cdots$, $v_n$ to 2. Since no two exterior vertices are adjacent in $T$, it is evident that $\sigma^*$ is a proper 3-coloring of $G$. So $\chi(G)\leq 3$.

{\textbf{Case 2.}}\quad $n$ is odd.

For the coloring $\sigma$, if there exists $j\in \{1, 2, \ldots, n\}$ such that $\sigma(v_j)\neq \sigma(v_{j+1})$, where $v_{n+1}=v_1$, let $v_{j+2}$, $v_{j+4}$, \ldots, $v_{j-1}$ be recolored with 2. It is easy to check that this yields a proper 3-coloring of $G$, and therefore $\chi(G)\leq 3$.

Otherwise, $\sigma(v_i)= \sigma(v_{i+1})$ for all $i\in \{1, 2, \ldots, n\}$, where $v_{n+1}=v_1$. Since $T$ is connected, $\sigma$ is uniquely determined up to swapping colors 0 and 1. Without loss of generality, let $\sigma(v_i)=1$ for all $i=1, 2, \ldots, n$. Because every vertex in $U$ is adjacent in $T$ to at least one vertex on $C_n$, all vertices in $U$ must be colored 0 under $\sigma$. Since $\sigma$ is a proper 2-coloring of $T$, $U$ is an independent set in $T$.

Because $E^{\prime}\neq \emptyset$, there exists at least one edge of $T$ joining two interior vertices. Since $U$ is an independent set in $T$, not all interior vertices can belong to $U$; thus there exists an interior vertex $w\in V(G)\setminus (U\cup V(C_n))$ (which has no neighbours on $C_n$). Let $uw\in E(T)$ be an edge incident to $w$, and let $T_u$ and $T_w$ be the two connected components of $T-uw$ containing $u$ and $w$, respectively. Since all leaves of $T$ lie on $C_n$, both $V(T_u)\cap V(C_n)$ and $V(T_w)\cap V(C_n)$ are non-empty, so there exists an edge $v_jv_{j+1}\in E(C_n)$ along the outer cycle with $v_j\in V(T_u)$ and $v_{j+1}\in V(T_w)$.

We now define a 3-coloring $\sigma^*$ of $G$ as follows:
\begin{itemize}
  \item For each $x\in V(T_u)$, let $\sigma^*(x)=\sigma(x)\in \{0, 1\}$ (in particular, $\sigma^*(v_j)=1$);
  \item For each $x\in V(T_w)\setminus \{w\}$, let $\sigma^*(x)=1-\sigma(x)\in \{0, 1\}$ (in particular, $\sigma^*(v_{j+1})=0$);
  \item Let $\sigma^*(w)=2$;
  \item Finally, recolor the $(n-1)/2$ exterior vertices $v_{j+2}, v_{j+4}, \ldots, v_{j-1}$ on $C_n\setminus \{v_j, v_{j+1}\}$ with color 2.
\end{itemize}
Since $w\notin U$ has no neighbours on $C_n$ and all neighbours of $w$ in $T$ have colors in $\{0, 1\}$, no neighbour of $w$ in $G$ has color 2. Moreover, all vertices in $U$ have colors in $\{0, 1\}$, and $\sigma^*(v_j)=1\neq 0=\sigma^*(v_{j+1})$. Thus $\sigma^*$ is a proper 3-coloring of $G$, and $\chi(G)=3$.
\end{proof}

\begin{lem}
If $G=T\cup C_n$ is a pseudo-Halin graph with $E^{\prime}= \emptyset$, then $$\chi(G)=\left\{\begin{array}{l}
4,\ \ \ \ \mbox{\rm if}\ G\in \mathcal{A};\\
3,\ \ \ \ \mbox{\rm otherwise}.\\
\end{array} \right.$$
\end{lem}

\begin{proof}
If $G=W_{p_1}$ is a wheel, then $\chi(W_{p_1})=4$ when $p_1$ is odd (i.e., $G\in \mathcal{A}$) and $\chi(W_{p_1})=3$ when $p_1$ is even (i.e., $G\notin \mathcal{A}$). For $r\geq 1$, by Lemma 2.3, let
$$H_1=R_r=W_{p_{r+1}}+_HW_{p_r},$$
$$H_2=R_{r-1}=H_1+_HW_{p_{r-1}},$$
$$\cdots$$
$$H_r=R_1=H_{r-1}+_HW_{p_1}.$$
Hence $H_1$, $H_2$, $\cdots$, $H_r$ are pseudo-Halin graphs with $E^{\prime}=\emptyset$.

First, $\chi(H_i)\geq3$, $i=1, 2, \ldots, r$. Next, we compute the chromatic number of $H_1$. Assume that the interior vertices of $H_1$ are $u_1$ and $u_2$, where $v_1, v_2, \ldots, v_{p_{r+1}}$ are adjacent to $u_1$, and $v_{p_{r+1}}, v_{p_{r+1}+1}, \ldots, v_{p_{r+1}+p_r-1}$ are adjacent to $u_2$.

\begin{figure}[htbp]
\centering
\includegraphics[height=4.5cm, width=0.24\textwidth]{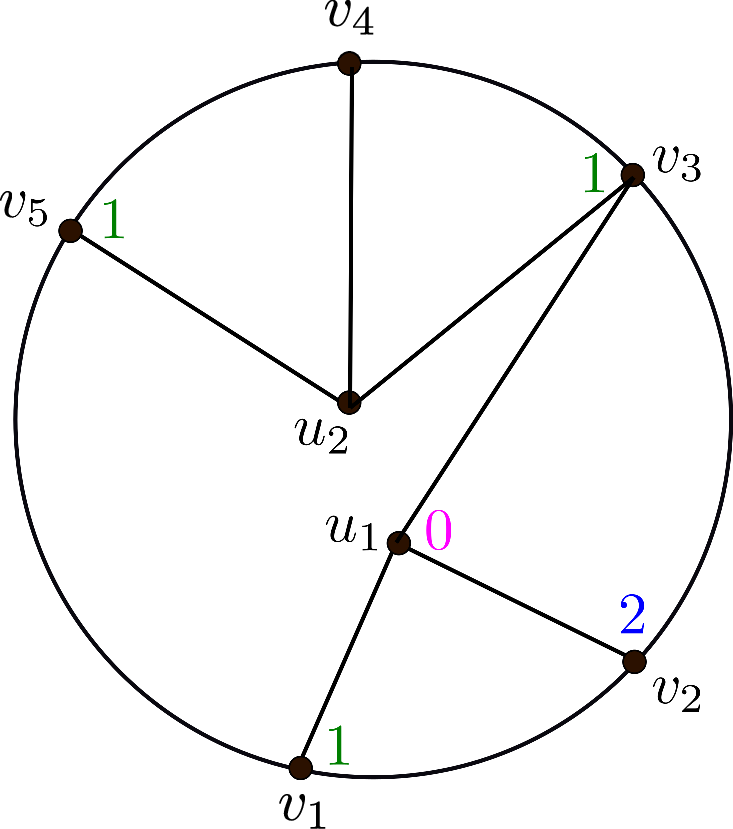}
\centerline{Fig.2  ~The pseudo-Halin graph $H_1$ with $p_{r+1}=p_r=3$.}
\end{figure}

When one of $p_{r+1}, p_r$ is odd and the other is even, $H_1$ has an even number $p_{r+1}+p_r-1$ of exterior vertices. It is easily seen that $\chi(H_1)=3$: let $u_1$, $u_2$ be colored 0; let $v_j$ be colored 1 if $j$ is odd, and 2 if $j$ is even, for $j=1, 2, \ldots, p_{r+1}+p_r-1$.

When $p_{r+1}, p_r$ are both even, let $v_1$, $v_3$, \ldots, $v_{p_{r+1}-1}$ be colored 0, $v_{p_{r+1}+1}$, $v_{p_{r+1}+3}$, \ldots, $v_{p_{r+1}+p_r-1}$ be colored 1, $v_2$, $v_4$, $\ldots$, $v_{p_{r+1}+p_r-2}$ be colored 2, $u_1$ be colored 1, and $u_2$ be colored 0. It is easy to check that the above coloring is a proper 3-coloring. Therefore $\chi(H_1)=3$.

When $p_{r+1}, p_r$ are both odd, without loss of generality, suppose $H_1$ admits a proper 3-coloring with $u_1$, $v_1$ colored $0$, $1$, respectively. Then the colors of $v_2$, \ldots, $v_{p_{r+1}}$ are uniquely forced, where the vertex $v_{p_{r+1}}$ is colored 1. The vertex $u_2$ is adjacent to $v_{p_{r+1}}$, so regardless of whether the color of $u_2$ is $0$ or $2$, the vertex $v_{p_{r+1}+p_r-1}$ must be colored $1$ (see Figure 2). Then the adjacent vertices $v_{p_{r+1}+p_r-1}$ and $v_1$ receive the same color, a contradiction. Hence $\chi(H_1)\geq4$. By the Four Color Theorem, $\chi(H_1)=4$.

Assume that the conclusion holds for $i\leq r-1$. Let $V(W_{p_1})=\{w_0, w_1, w_2, \ldots, w_{p_1}\}$, where $w_0$ is the center vertex. For the pseudo-Halin graph $H_{r-1}=C_{n^{\prime}}\cup T^{\prime}$, let the outer cycle of $H_{r-1}$ be $v_1v_2\ldots v_{n^{\prime}}v_1$. Without loss of generality, $H_r=H_{r-1}(v_jv_{j+1})+_HW_{p_1}(w_1w_{p_1})$ is obtained by identifying the vertices $v_j$, $w_1$ as $v_H(v_jw_1)$ for some $j\in \{1, 2, \cdots, n^{\prime}\}$ and adding the edge $v_{j+1}w_{p_1}$. We write $\sigma^*$ for a coloring of $H_r$.

If $H_{r-1}\notin \mathcal{A}$, by induction, $\chi(H_{r-1})=3$. It may be assumed that $H_{r-1}$ has a proper 3-coloring $\sigma$ with $\sigma(v_j)=0, \sigma(v_{j+1})=1$. Let $\sigma^*(v_H(v_jw_1))=0$ and $\sigma^*(w_0)=1$; then for $l=2, 3, \ldots, p_1$, let $\sigma^*(w_l)=2$ if $l$ is even and $\sigma^*(w_l)=0$ if $l$ is odd, while $\sigma^*(v)=\sigma(v)$ for each $v\in V(H_{r-1})\setminus \{v_j\}$. Since $\sigma^*(w_{p_1})\in \{0, 2\}$ and $\sigma^*(v_{j+1})=1$, it is evident that $\sigma^*$ is proper, so $\chi (H_{r})=3$.

If $H_{r-1}\in \mathcal{A}$, by induction, $\chi(H_{r-1})=4$, and all $p_2, \ldots, p_{r+1}$ are odd, which implies that $H_{r-1}$ has an odd number $n^{\prime}=\sum_{i=2}^{r+1}(p_i-1)+1$ of exterior vertices. When $p_1$ is even, $H_r=T\cup C_n$ has $n=n^{\prime}+p_1-1$ exterior vertices, which is even. Let $\sigma^*(u)=0$ for all $u\in V(T)\setminus V(C_n)$, and let $\sigma^*(v_i)=1$ when $i$ is odd and $\sigma^*(v_i)=2$ when $i$ is even, for $v_i\in V(C_n)$ $(i=1, 2, \ldots, n)$. Since $E^{\prime}=\emptyset$, $\sigma^*$ is proper and $\chi(H_r)=3$. When $p_1$ is odd, assume that $\sigma^*$ is a proper 3-coloring of $H_r$; since $p_1$ is odd, the path $v_H(v_jw_1)w_2\ldots w_{p_1}$ in the wheel $W_{p_1}$ forces $\sigma^*(v_H(v_jw_1))=\sigma^*(w_{p_1})\neq \sigma^*(v_{j+1})$. Then a proper 3-coloring $\sigma$ of $H_{r-1}$ is obtained by restricting $\sigma^*$ to $V(H_{r-1})$ with $\sigma(v_j)=\sigma^*(v_H(v_jw_1))$, contradicting $\chi(H_{r-1})=4$. Hence $\chi(H_r)\geq 4$, and by the Four Color Theorem, $\chi(H_r)=4$.

This completes the proof.
\end{proof}

\noindent
\textbf{Proof of Theorem 1.}~
Combining Lemma 3.1 (for the case $E^{\prime}\neq \emptyset$) and Lemma 3.2 (for the case $E^{\prime}=\emptyset$), we immediately obtain that $\chi(G)=4$ if $G\in \mathcal{A}$ and $\chi(G)=3$ otherwise. This completes the proof of Theorem 1. \hfill $\Box$

\section{Proof of Theorem 2}
In this section, we divide the proof of Theorem 2 into two parts according to whether $G\in \mathcal{A}$ or $G\notin \mathcal{A}$.

\subsection*{Part 1: The case $G\in \mathcal{A}$}

For any subgraph $H$ of a pseudo-Halin graph $G=T\cup C_n$, we have
$$|E(H)|=|E(H)\cap E(T)|+|E(H)\cap E(C_n)|\leq (|V(H)|-1)+|V(H)|=2|V(H)|-1,$$
which implies $\delta(H)\leq \frac{2|E(H)|}{|V(H)|}<4$, and thus $\delta(H)\leq 3$. Since $\delta(G)\geq 3$, the degeneracy of any pseudo-Halin graph $G$ is 3. By Lemma 2.1, $AT(G)\leq 4$ holds for every pseudo-Halin graph $G$. When $G\in \mathcal{A}$, Theorem 1 gives $\chi(G)=4$, and hence $AT(G)\geq \chi(G)=4$. Consequently, $AT(G)=4$ for every pseudo-Halin graph $G\in \mathcal{A}$.

\subsection*{Part 2: The case $G\notin \mathcal{A}$}

When a pseudo-Halin graph $G=T\cup C_n$ is not in $\mathcal{A}$, Theorem 1 gives $\chi(G)=3$, and hence $AT(G)\geq \chi(G)=3$. To prove $AT(G)=3$, it remains to establish that $AT(G)\leq 3$. We first recall the Alon-Tarsi number of wheel graphs.

\begin{lem}{\cite{ZG2021}}
If $G$ is a wheel graph $W_p$ of order $p+1$, then
$$AT(W_p)=\left\{\begin{array}{l}
4,\ \ \ \mbox{\rm if } p\ \mbox{\rm is odd};\\
3,\ \ \ \mbox{\rm if } p\ \mbox{\rm is even}.\\
\end{array} \right.$$
\end{lem}

If $G\notin \mathcal{A}$ is a wheel, then $G=W_p$ with $p$ even, so $AT(G)=3$ by Lemma 4.1. When $G$ is neither in $\mathcal{A}$ nor a wheel, let $P=u_1u_2\cdots u_k$ be a longest path in $T$ and take $w\in \{u_2, u_{k-1}\}$ with $d_G(w)=m+1$ $(m\geq 2)$. By Lemma 2.2, $w$ satisfies either Condition (1) or Condition (2) of Lemma 2.2. We prove $AT(G)\leq 3$ by induction on $|V(G)|$ through the following two lemmas (which proceed by simultaneous induction on $|V(G)|$, with even wheels in Lemma 4.1 serving as base cases whenever a reduction yields a wheel).

\begin{lem}
If a pseudo-Halin graph $G=T\cup C_n$ is neither in $\mathcal{A}$ nor a wheel and satisfies Condition $(1)$ of Lemma 2.2, then $AT(G)\leq3$.
\end{lem}

\begin{proof}
We proceed by induction on the number of vertices $|V(G)|$. Since $G$ is not a wheel, $G$ has at least two interior vertices. For an illustrative small graph $G$ (e.g., when $m=4$, see Figure 3 (a)), let $D$ be an orientation of $G$ obtained by orienting the outer cycle $C_8$ clockwise and every $xv_5$-path in $T$ from $x$ to $v_5$ for $x\in V(G)\setminus\{v_5\}$. Denote the directed path from $v_j$ to $v_5$ by $v_j\rightsquigarrow v_5$. Since $d_D^+(v_5)=1$, every directed cycle in $D$ other than $C_8$ must pass through $v_5$ and consists of $v_j\rightsquigarrow v_5$ in $T$ and $v_5\rightsquigarrow v_j$ in $C_8$ for $j\in \{1, 2, 3, 4, 6, 7, 8\}$. Consequently, no two directed cycles in $D$ are edge-disjoint, and $\mathcal{E}(D)$ consists of the empty subdigraph, $C_8$, and the 7 directed cycles above, so $|\mathcal{E}(D)|=9$ is odd. Thus ${\rm diff}(D)\neq 0$, and $D$ is an $AT$-orientation of $G$ with maximum outdegree 2 (see Figure 3 (b)), giving $AT(G)\leq 3$.

\begin{figure}[htbp]
\centering
\includegraphics[height=5cm, width=0.56\textwidth]{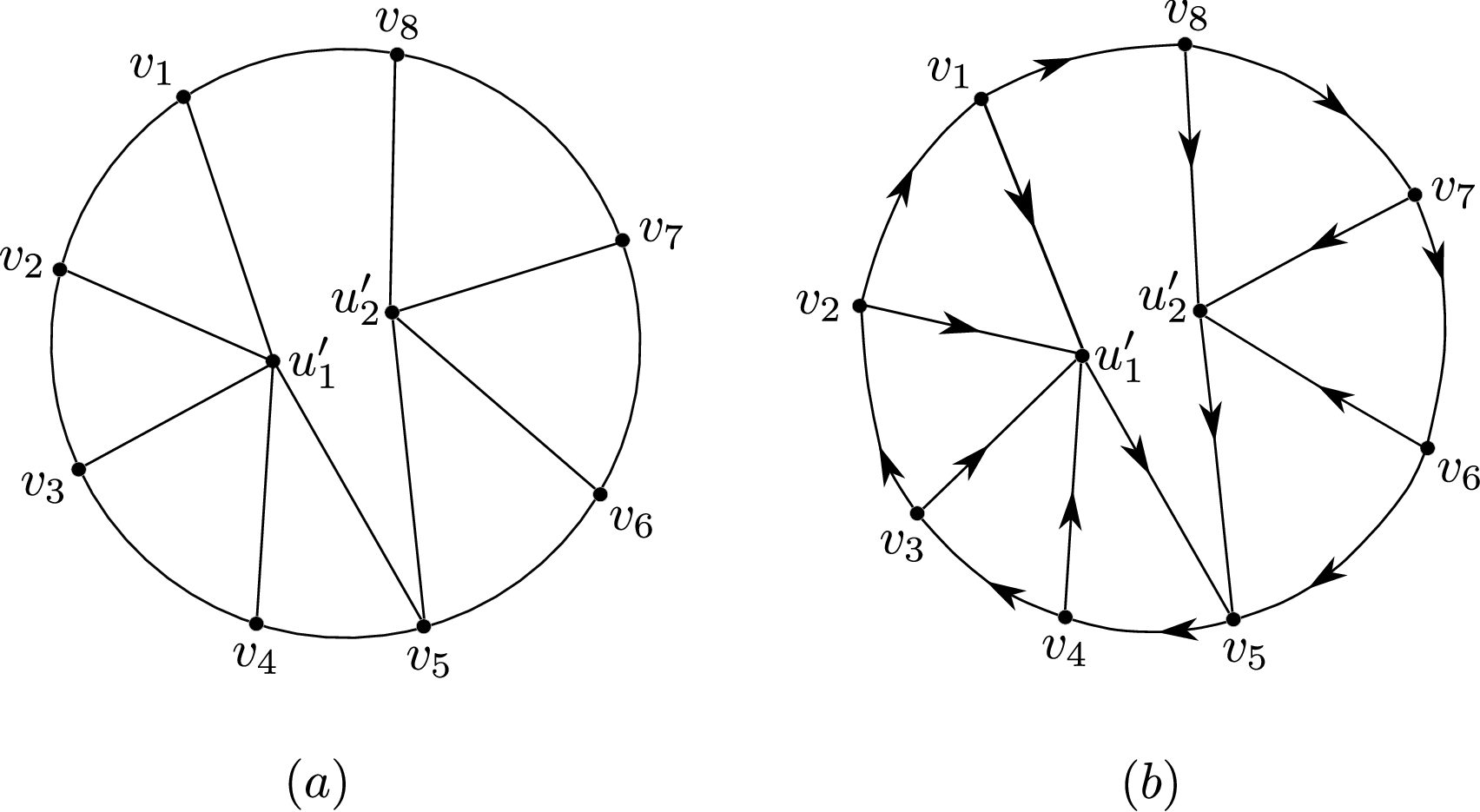}
\centerline{Fig.3  ~$(a)$ The pseudo-Halin graph $G$ with ten vertices.   $(b)$ The orientation of $G$.}
\end{figure}

Now assume that $AT(G^{\prime})\leq 3$ holds for every pseudo-Halin graph $G^{\prime}\notin \mathcal{A}$ with $|V(G^{\prime})|<|V(G)|$. Let $N(w)=\{v_1, v_2, \cdots, v_m, y\}$ $(m\geq 2)$ be the set of neighbours of $w$ in $G$, where $v_1, \cdots, v_m\in R(f_0)$, $y\in IR(f_0)$, and $x, v_1, v_2, \cdots, v_m, y$ lie on $C_n$ in counterclockwise order.

\begin{figure}[htbp]
\centering
\includegraphics[height=7cm, width=0.65\textwidth]{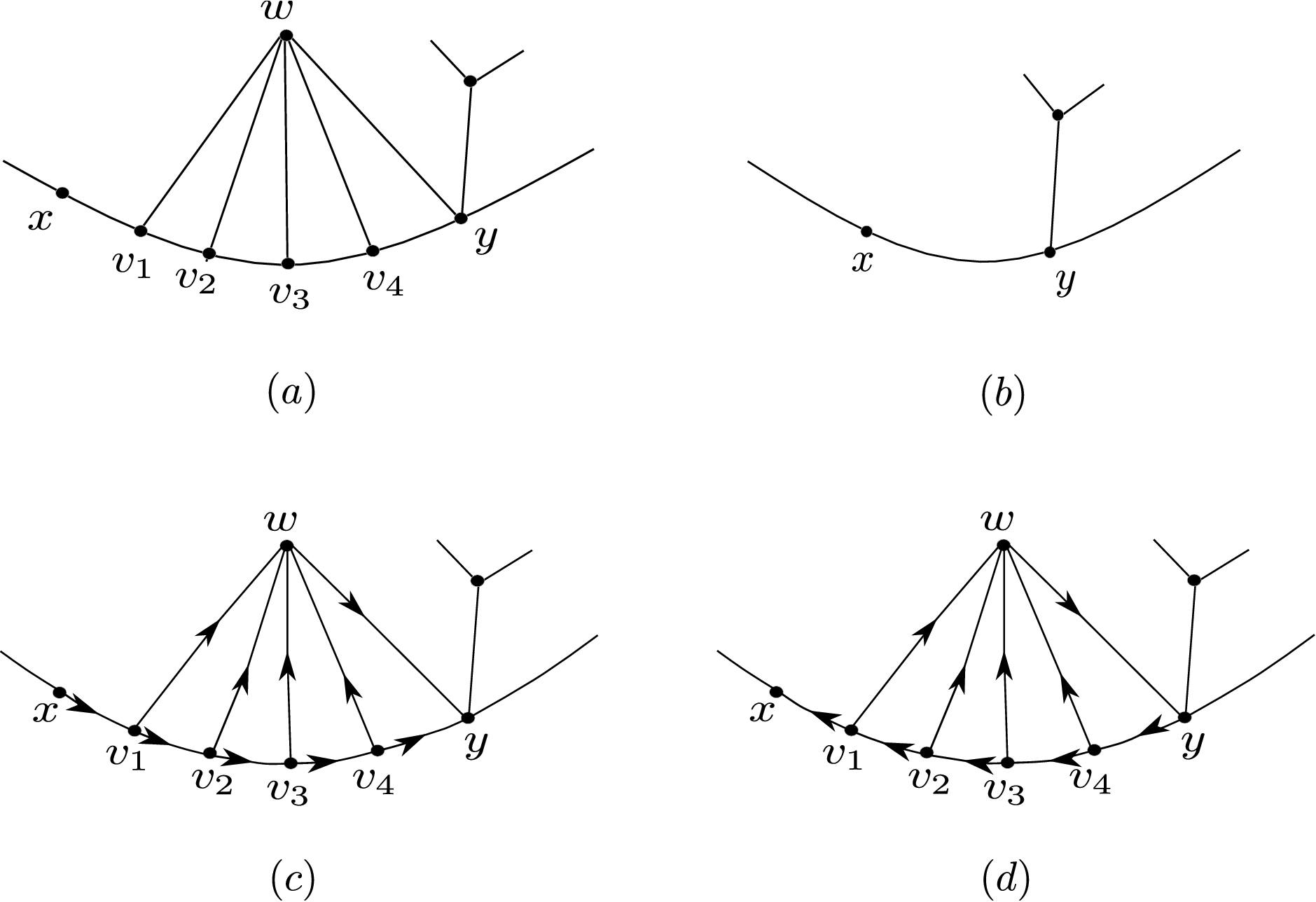}
\centerline{Fig.4 (for $m=4$)  ~$(a)$ The graph $G$.   $(b)$ The graph $G_1$.}
\centerline{$(c)$ The orientation $D$ when the edge $xy$ is oriented from $x$ to $y$ in $D_1$.}
\centerline{$(d)$ The orientation $D$ when the edge $xy$ is oriented from $y$ to $x$ in $D_1$.}
\end{figure}

{\textbf{Case 1.}}\quad $m\geq 2$ is even.

Let $G_1=G-\{w, v_i \mid i=1, 2, \cdots ,m\}+\{xy\}$. By Lemma 2.2(1), $G_1$ is a pseudo-Halin graph. Since $m$ is even, $d_G(w)=m+1$ is odd; because $G\notin \mathcal{A}$ and all interior vertices of $G_1$ have the same degrees in $G_1$ as in $G$ with $E^{\prime}(G_1)=E^{\prime}(G)$, we have $G_1\notin \mathcal{A}$. By the induction hypothesis (or Lemma 4.1 if $G_1$ is an even wheel), $G_1$ has an $AT$-orientation $D_1$ with $d_{D_1}^+(v)\leq 2$ for each $v\in V(G_1)$. Based on $D_1$, we construct an orientation $D$ of $G$.

If $xy$ is oriented from $x$ to $y$ in $D_1$, let $D$ be obtained from $D_1$ by deleting $(x, y)$ and adding the arcs $(v_i, w)$ for $i=1, 2, \cdots, m$, $(w, y)$, $(x, v_1)$, $(v_m, y)$, and $(v_j, v_{j+1})$ for $j=1,2, \cdots, m-1$ (see Figure 4 (c)). In $D$, no directed cycle lies entirely inside $\{w, v_1, \cdots, v_m, y\}$. For any Eulerian subdigraph of $D_1$ containing $(x, y)$, replacing the single arc $(x, y)$ by the path $xv_1v_2\cdots v_my$ of odd length $m+1$ preserves parity, whereas the $m$ directed paths $xv_1\cdots v_iwy$ $(i=1, \cdots, m)$ from $x$ to $y$ have lengths $3, 4, \cdots, m+2$, consisting of $\frac{m}{2}$ odd paths and $\frac{m}{2}$ even paths which cancel out in ${\rm diff}(D)$. Thus ${\rm diff}(D)=|\mathcal{E}_e(D)|-|\mathcal{E}_o(D)|=|\mathcal{E}_e(D_1)|-|\mathcal{E}_o(D_1)|\neq 0$.

If $xy$ is oriented from $y$ to $x$ in $D_1$, let $D$ be obtained from $D_1$ by deleting $(y, x)$ and adding the arcs $(v_i, w)$ for $i=1, 2, \cdots, m$, $(w, y)$, $(v_1, x)$, $(y, v_m)$, and $(v_{j+1}, v_j)$ for $j=1,2, \cdots, m-1$ (see Figure 4 (d)). Here, replacing $(y, x)$ by the path $yv_mv_{m-1}\cdots v_1x$ of odd length $m+1$ preserves parity, while the $m$ directed cycles $yv_m\cdots v_iwy$ $(i=1, \cdots, m)$ containing $(w, y)$ and $(y, v_m)$ consist of $\frac{m}{2}$ odd cycles and $\frac{m}{2}$ even cycles which cancel out. Thus ${\rm diff}(D)={\rm diff}(D_1)\neq0$.

In both cases, $d_{D}^+(v)\leq 2$ for all $v\in V(G)$, so $D$ is an $AT$-orientation of $G$ and $AT(G)\leq 3$.

\begin{figure}[htbp]
\centering
\includegraphics[height=17.6cm, width=0.65\textwidth]{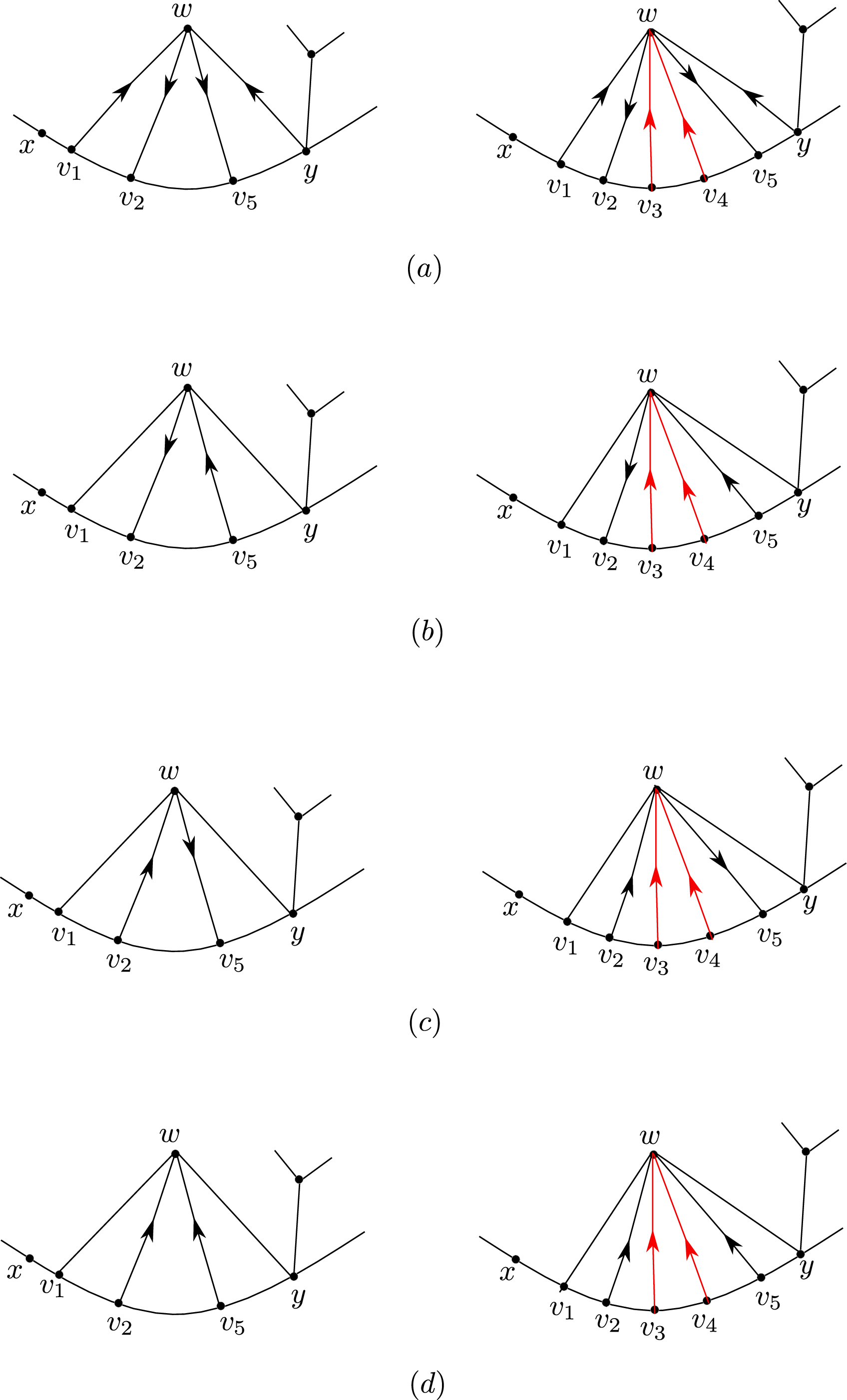}
\centerline{Fig.5 (for $m=5$)  ~$(a)$-$(d)$ The orientations of $G_2$ and $G$ in Case 2.}
\end{figure}

{\textbf{Case 2.}}\quad $m\geq 3$ is odd.

{\textbf{Subcase 2.1.}}\quad $m=3$.

Let $G_1=G-\{w, v_1, v_2, v_3\}+\{xy\}$, which is a pseudo-Halin graph by Lemma 2.2(1).
\begin{itemize}
  \item If $G_1\notin \mathcal{A}$, then by the induction hypothesis (or Lemma 4.1), $G_1$ admits an $AT$-orientation $D_1$ with $d_{D_1}^+(v)\leq 2$. Let $\mathcal{E}^0(D_1)$ and $\mathcal{E}^1(D_1)$ be the sets of Eulerian subdigraphs of $D_1$ not containing $xy$ and containing $xy$, respectively, and let $\Delta_0=|\mathcal{E}^0_e(D_1)|-|\mathcal{E}^0_o(D_1)|$ and $\Delta_1=|\mathcal{E}^1_e(D_1)|-|\mathcal{E}^1_o(D_1)|$. Since ${\rm diff}(D_1)=\Delta_0+\Delta_1\neq 0$, $\Delta_0$ and $\Delta_1$ cannot both be 0.
  \begin{enumerate}
    \item[(i)] If $\Delta_0\neq 0$, let $D$ be obtained from $D_1-\{xy\}$ by orienting the fan edges as $(v_1, x)$, $(v_1, v_2)$, $(v_2, v_3)$, $(v_3, y)$, $(w, v_1)$, $(v_2, w)$, $(v_3, w)$, and $(w, y)$. Then $d_D^+(v)\leq 2$ for all $v\in V(G)$. Since $(v_1, x)$, $(v_3, y)$, and $(w, y)$ are all oriented from $\{w, v_1, v_2, v_3\}$ toward $\{x, y\}$, no directed cycle in $D$ can cross between $\{w, v_1, v_2, v_3\}$ and $V(G_1)$. Inside $D[\{w, v_1, v_2, v_3\}]$, the odd 3-cycle $wv_1v_2w$ and the even 4-cycle $wv_1v_2v_3w$ cancel out, leaving only the empty subdigraph. Hence ${\rm diff}(D)=\Delta_0\neq 0$.
    \item[(ii)] If $\Delta_0=0$, then $\Delta_1\neq 0$. When $(x, y)\in A(D_1)$, delete $(x, y)$ and orient the fan edges as $(x, v_1)$, $(v_1, v_2)$, $(v_2, v_3)$, $(v_3, y)$, $(v_1, w)$, $(v_2, w)$, $(w, v_3)$, and $(w, y)$; this fan is acyclic and contains 3 even paths and 2 odd paths from $x$ to $y$, giving ${\rm diff}(D)=\Delta_0-\Delta_1=-\Delta_1\neq 0$. Symmetrically, when $(y, x)\in A(D_1)$, delete $(y, x)$ and orient the fan edges as $(y, v_3)$, $(v_3, v_2)$, $(v_2, v_1)$, $(v_1, x)$, $(v_1, w)$, $(v_2, w)$, $(w, v_3)$, and $(w, y)$; the 4 internal directed cycles (two of length 4, one of length 3, one of length 5) all share $(v_3, v_2)$ and cancel out, while the path $yv_3v_2v_1x$ has even length 4 and blocks all internal cycles, again giving ${\rm diff}(D)=\Delta_0-\Delta_1=-\Delta_1\neq 0$.
  \end{enumerate}
  \item If $G_1\in \mathcal{A}$, then by Lemma 2.3, $G_1$ is an iterated Haj\'{o}s sum of wheels. Deleting the exterior cycle edge $xy$ makes $G_1-\{xy\}=G-\{w, v_1, v_2, v_3\}$ a 2-degenerate graph (as vertices of degree at most 2 can be peeled successively starting from $x$ and $y$ across each wheel). Thus $G_1-\{xy\}$ admits an acyclic orientation $D^{\prime}_1$ with maximum outdegree at most 2 and ${\rm diff}(D^{\prime}_1)=1$. Combining $D^{\prime}_1$ with the fan orientation in (i) above yields an orientation $D$ of $G$ with $d_D^+(v)\leq 2$ and ${\rm diff}(D)=1\cdot {\rm diff}(D^{\prime}_1)=1\neq 0$.
\end{itemize}

{\textbf{Subcase 2.2.}}\quad $m\geq 5$ is odd.

Let $G_2=G-\{v_3, \cdots, v_{m-1}\}+\{v_2v_m\}$. By Lemma 2.2(1), $G_2$ is a pseudo-Halin graph in which $w$ has 3 regular exterior neighbours $\{v_1, v_2, v_m\}$ (so $G_2$ falls into Subcase 2.1 with $m=3$, and since $d_{G_2}(w)=4$ is even, $G_2\notin \mathcal{A}$). By Subcase 2.1, $G_2$ has an $AT$-orientation $D_2$ with maximum outdegree 2. Let $D$ be the orientation of $G$ that orients $v_jw$ as $a_j=(v_j, w)$ for $j=3, \cdots, m-1$, gives $v_jv_{j+1}$ $(j=2, \cdots, m-1)$ the same orientation as $v_2v_m$ in $D_2$, and orients all other edges as in $D_2$ (see Figure 5). Clearly, $d_{D}^+(v)\leq 2$ for every $v\in V(G)$.

To see that ${\rm diff}(D)\neq 0$, note first that any Eulerian subdigraph of $D$ not using $\{a_3, \cdots, a_{m-1}\}$ corresponds to an Eulerian subdigraph of $D_2$ in which the single arc on $v_2v_m$ is replaced by the directed path $v_2v_3\cdots v_m$ of odd length $m-2$, which preserves edge-count parity. Second, since each vertex $v_i$ $(3\leq i\leq m-1)$ has indegree 1 and outdegree 2 in $D$, any Eulerian subdigraph $H\in \mathcal{E}(D)$ can contain at most one arc $a_i=(v_i, w)$ for $3\leq i\leq m-1$, and the sets $\mathcal{E}(D, a_i)$ of Eulerian subdigraphs containing $a_i$ are pairwise disjoint. Moreover, in any $H\in \mathcal{E}(D, a_i)$, the arcs of $H$ incident to $\{v_3, \cdots, v_{m-1}\}$ consist precisely of the directed rim path from $v_2$ (or $v_m$) to $v_i$ followed by $a_i$. Since $m-3$ is even, pairing $i=2k+1$ with $i+1=2k+2$ for $k=1, \cdots, \frac{m-3}{2}$ changes the path length by $\pm 1$ and thus establishes a bijection between $\bigcup_{i=3}^{m-1}\mathcal{E}_e(D, a_i)$ and $\bigcup_{i=3}^{m-1}\mathcal{E}_o(D, a_i)$. Therefore,
$${\rm diff}(D)=\left(|\mathcal{E}_e(D_2)|+\sum\limits_{i=3}^{m-1}|\mathcal{E}_e(D, a_i)|\right)-\left(|\mathcal{E}_o(D_2)|+\sum\limits_{i=3}^{m-1}|\mathcal{E}_o(D, a_i)|\right)={\rm diff}(D_2)\neq 0.$$
This completes the proof of Lemma 4.2.
\end{proof}

\begin{lem}
If a pseudo-Halin graph $G=T\cup C_n$ is neither in $\mathcal{A}$ nor a wheel and satisfies Condition $(2)$ of Lemma 2.2, then $AT(G)\leq3$.
\end{lem}

\begin{proof}
We proceed by induction on $|V(G)|$. For an illustrative small graph $G$ satisfying Condition (2) of Lemma 2.2 (e.g., the 8-vertex graph in Figure 6 (a)), let $D$ be the orientation of $G$ obtained by orienting the outer cycle $C_6$ clockwise and every $xv_6$-path in $T$ from $x$ to $v_6$ for $x\in V(G)\setminus \{v_6\}$ (see Figure 6 (b)). Since $d_D^+(v_6)=1$, every directed cycle in $D$ except $C_6$ consists of $v_j\rightsquigarrow v_6$ in $T$ and $v_6\rightsquigarrow v_j$ along $C_6$ for $j=1, 2, 3, 4, 5$. Together with the empty subdigraph and $C_6$, we have $|\mathcal{E}(D)|=7$ (odd), so ${\rm diff}(D)\neq 0$ and $d_D^+(v)\leq 2$ for all $v\in V(G)$, giving $AT(G)\leq 3$.

\begin{figure}[htbp]
\centering
\includegraphics[height=4.8cm, width=0.58\textwidth]{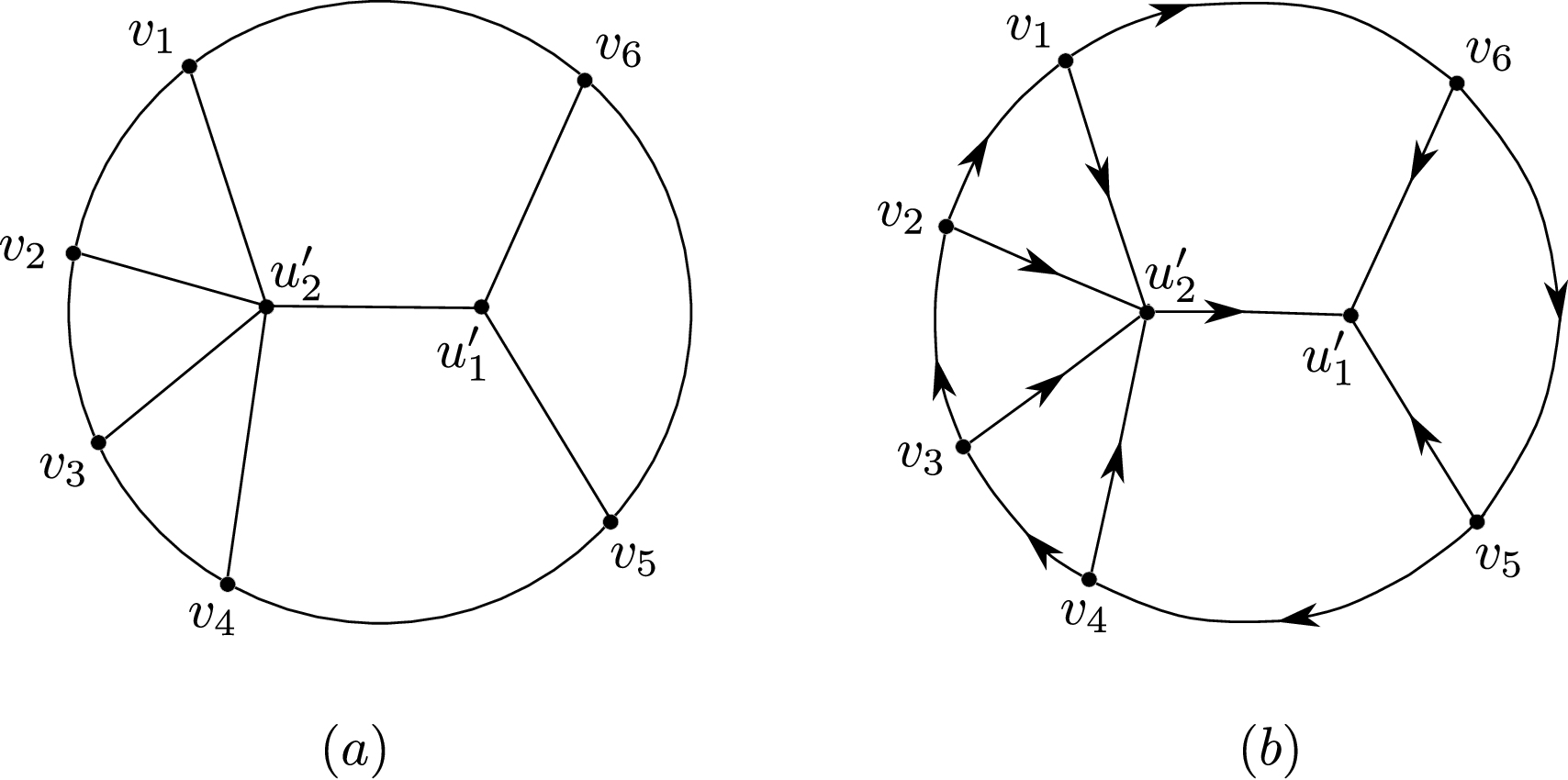}
\centerline{Fig.6  ~$(a)$ The pseudo-Halin graph $G$ with eight vertices.   $(b)$ The orientation of $G$.}
\end{figure}

Now assume that $AT(G^{\prime})\leq 3$ holds for every pseudo-Halin graph $G^{\prime}\notin \mathcal{A}$ with $|V(G^{\prime})|<|V(G)|$. Let $N(w)=\{u, v_1, v_2, \cdots, v_m\}$ $(m\geq 2)$, where $u$ is an interior vertex, $v_1, \cdots, v_m\in R(f_0)$, and $x, v_1, v_2, \cdots, v_m, y$ lie on $C_n$ in counterclockwise order.

{\textbf{Case 1.}}\quad $m\geq 3$ is odd.

{\textbf{Subcase 1.1.}}\quad $m=3$.

By Lemma 2.2(2), $G_2^1=G-\{v_1, v_2, v_3\}+\{xw, yw\}$ is a pseudo-Halin graph in which $w$ is a regular exterior vertex of degree 3 incident to $xw, yw, uw$.
\begin{itemize}
  \item If $G_2^1\in \mathcal{A}$, then deleting the degree-3 regular exterior vertex $w$ from $G_2^1\in \mathcal{A}$ yields $G-\{w, v_1, v_2, v_3\}=G_2^1-\{w\}$, which is 2-degenerate and therefore admits an acyclic orientation $D^{\prime}$ with maximum outdegree at most 2 and ${\rm diff}(D^{\prime})=1$. Extending $D^{\prime}$ to $G$ by orienting the fan edges as $(v_1, x)$, $(v_1, v_2)$, $(v_2, v_3)$, $(v_3, y)$, $(w, v_1)$, $(v_2, w)$, $(v_3, w)$, and $(w, u)$ ensures that all three cut-arcs $(v_1, x)$, $(v_3, y)$, $(w, u)$ leave $\{w, v_1, v_2, v_3\}$, while the internal 3-cycle $wv_1v_2w$ and 4-cycle $wv_1v_2v_3w$ cancel out. Thus ${\rm diff}(D)=1\cdot {\rm diff}(D^{\prime})=1\neq 0$.
  \item If $G_2^1\notin \mathcal{A}$, then by the induction hypothesis (or Lemma 4.1), $G_2^1$ admits an $AT$-orientation $D^{\prime}$ with $d_{D^{\prime}}^+(v)\leq 2$. In $D^{\prime}$, the Eulerian subdigraphs partition into those not visiting $w$ (contributing $\Delta_0$) and those passing through $w$ via the at most two directed paths $f_1, f_2$ of length 2 (contributing $\Delta_{f_1}$ and $\Delta_{f_2}$), so that ${\rm diff}(D^{\prime})=\Delta_0+\Delta_{f_1}+\Delta_{f_2}\neq 0$. If $\Delta_0\neq 0$, the fan orientation in the previous bullet gives ${\rm diff}(D)=\Delta_0\neq 0$. If $\Delta_0=0$, then $\Delta_{f_1}+\Delta_{f_2}\neq 0$, and replacing $w$ by the complementary fan orientation on $\{w, v_1, v_2, v_3\}$ (for instance, when $(x, w)$, $(y, w)$, $(w, u)\in A(D^{\prime})$, orienting $(x, v_1)$, $(v_1, v_2)$, $(v_3, v_2)$, $(y, v_3)$, $(v_1, w)$, $(v_2, w)$, $(v_3, w)$, $(w, u)$, which has 2 even and 1 odd paths from $x$ to $u$ and from $y$ to $u$) yields ${\rm diff}(D)=\Delta_{f_1}+\Delta_{f_2}\neq 0$.
\end{itemize}

{\textbf{Subcase 1.2.}}\quad $m\geq 5$ is odd.

Let $G_3=G-\{v_3, \cdots, v_{m-1}\}+\{v_2v_m\}$. By Lemma 2.2(2), $G_3$ is a pseudo-Halin graph in which $w$ has $m=3$ regular exterior neighbours $\{v_1, v_2, v_m\}$, and $E^{\prime}(G_3)=E^{\prime}(G)\neq \emptyset$ so $G_3\notin \mathcal{A}$. By Subcase 1.1, $G_3$ has an $AT$-orientation $D_3$ with maximum outdegree 2. Construct the orientation $D$ of $G$ from $D_3$ by orienting $v_jw$ as $a_j=(v_j, w)$ for $j=3, \cdots, m-1$, giving $v_jv_{j+1}$ $(j=2, \cdots, m-1)$ the same direction as $v_2v_m$ in $D_3$, and keeping all other arcs as in $D_3$ (see Figure 7).

\begin{figure}[htbp]
\centering
\includegraphics[height=3.6cm, width=0.7\textwidth]{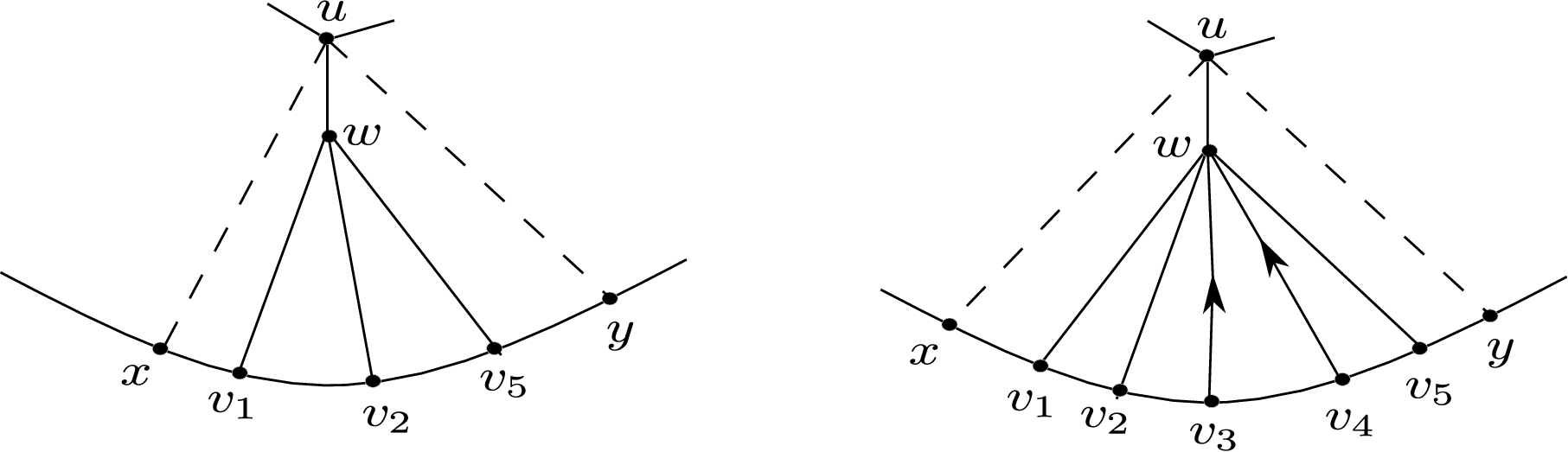}
\centerline{Fig.7 (for $m=5$) ~The graphs $G_3$ and $G$.}
\end{figure}

By the exact same pairing bijection as in Subcase 2.2 of Lemma 4.2, the path $v_2v_3\cdots v_m$ has odd length $m-2$ (preserving the parity of Eulerian subdigraphs of $D_3$), and pairing $a_{2k+1}$ with $a_{2k+2}$ for $k=1, \cdots, \frac{m-3}{2}$ cancels all Eulerian subdigraphs using $\{a_3, \cdots, a_{m-1}\}$. Hence ${\rm diff}(D)={\rm diff}(D_3)\neq 0$.

{\textbf{Case 2.}}\quad $m\geq 2$ is even.

{\textbf{Subcase 2.1.}}\quad $m=2$.

Consider $G_0=G-\{w, v_1, v_2\}+\{xy\}$ and $G_2^1=G-\{v_1, v_2\}+\{xw, yw\}$. For any orientation $D_0$ of $G_0$ with maximum outdegree at most 2, if we orient $(v_1, w)$, $(v_2, w)$, $(w, u)$, and give the path $xv_1v_2y$ the same direction as $xy$ in $D_0$, then $(w, u)$ enters $u$ (adding 0 to $d^+(u)$), the two paths of lengths 3 and 4 from $x$ (or $y$) through $w$ to $u$ cancel out, and the path $xv_1v_2y$ has odd length 3, so ${\rm diff}(D)={\rm diff}(D_0)$.
\begin{itemize}
  \item If $d_G(u)\geq 4$, then $d_{G_0}(u)=d_G(u)-1\geq 3$ and $G_0$ is a pseudo-Halin graph. If $G_0\notin \mathcal{A}$, then by the induction hypothesis (or Lemma 4.1), $G_0$ has an $AT$-orientation $D_0$ and ${\rm diff}(D)={\rm diff}(D_0)\neq 0$. If $G_0\in \mathcal{A}$, then $E^{\prime}(G_0)=\emptyset$ and $d_{G_0}(u)$ is odd, so $d_{G_2^1}(u)=d_G(u)$ is even and $G_2^1\notin \mathcal{A}$ has an even outer cycle. Orienting $G_2^1$ so that $(w, u)$ is not in any directed cycle ($\Delta_{xu}=0, \Delta_0=1$) and using the fan orientation $(x, v_1)$, $(v_1, v_2)$, $(v_2, y)$, $(v_1, w)$, $(w, v_2)$, $(w, u)$ (where the two $x\to y$ paths of lengths 3 and 4 cancel) gives ${\rm diff}(D)=\Delta_0-\Delta_{xu}=1\neq 0$.
  \item If $d_G(u)=3$, then by Lemma 2.2(2), $G_2^1$ is a pseudo-Halin graph with $G_2^1\notin \mathcal{A}$ (or if $G_2^1\in \mathcal{A}$, $G_2^1=W_3$ and $G$ is the 6-vertex prism/Halin graph which has $AT(G)=3$ by $\cite{ZG2021}$). For an $AT$-orientation $D^{\prime}$ of $G_2^1$ with ${\rm diff}(D^{\prime})=\Delta_0+\Delta_{f_1}+\Delta_{f_2}\neq 0$ (where at least one of $\Delta_0, \Delta_{f_1}, \Delta_{f_2}$ differs from the others), the three fan orientations $D^{(1)}$ (with $(v_1, w), (v_2, w)$, giving $\Delta_0-\Delta_{f_1}$), $D^{(2)}$ (with $(v_1, w), (w, v_2)$, giving $\Delta_0-\Delta_{f_2}$), and $D^{(3)}$ (with $(w, v_1), (v_2, w)$, giving $\Delta_{f_2}-\Delta_{f_1}$) all satisfy $d^+(v)\leq 2$ and at least one has ${\rm diff}(D^{(i)})\neq 0$.
\end{itemize}

{\textbf{Subcase 2.2.}}\quad $m\geq 4$ is even.

Let $G_4=G-\{v_2, \cdots ,v_{m-1}\}+\{v_1v_m\}$. By Lemma 2.2(2), $G_4$ is a pseudo-Halin graph in which $w$ has $m=2$ regular exterior neighbours $\{v_1, v_m\}$ and $E^{\prime}(G_4)=E^{\prime}(G)\neq \emptyset$ (so $G_4\notin \mathcal{A}$). By Subcase 2.1, $G_4$ has an $AT$-orientation $D_4$ satisfying $d_{D_4}^+(v)\leq 2$ for each $v\in V(G_4)$. Extend $D_4$ to $G$ so that the orientation of $v_jv_{j+1}$ $(j=1, 2, \cdots, m-1)$ is the same as that of $v_1v_m$ in $D_4$, the edge $v_iw$ $(i=2, \cdots, m-1)$ is oriented as $a_i=(v_i, w)$, and all other edges are oriented as in $D_4$ (see Figures 8 and 9).

\begin{figure}[htbp]
\centering
\includegraphics[height=7.8cm, width=0.95\textwidth]{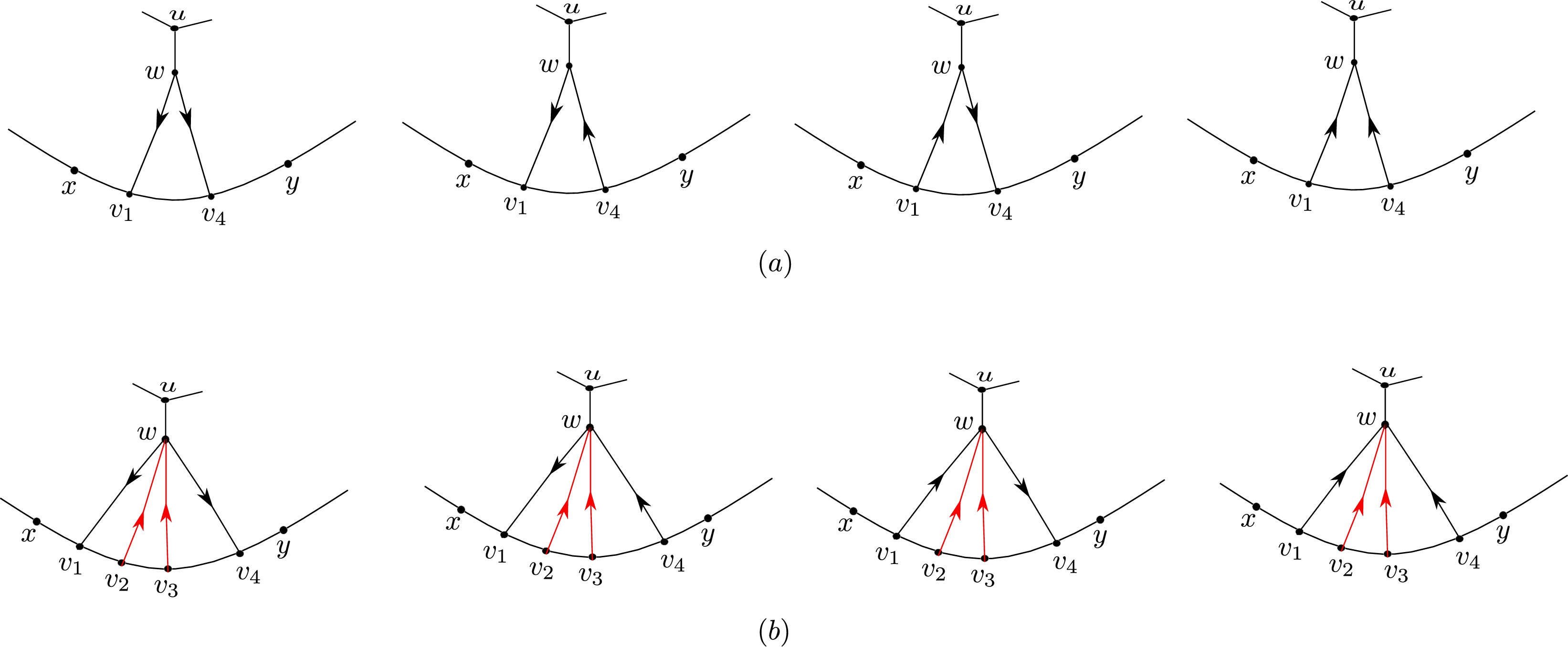}
\centerline{Fig.8 (for $m=4$)  ~$(a)$ The graph $G_4$.   $(b)$ The graph $G$.}
\end{figure}

\begin{figure}[htbp]
\centering
\includegraphics[height=3.4cm, width=0.7\textwidth]{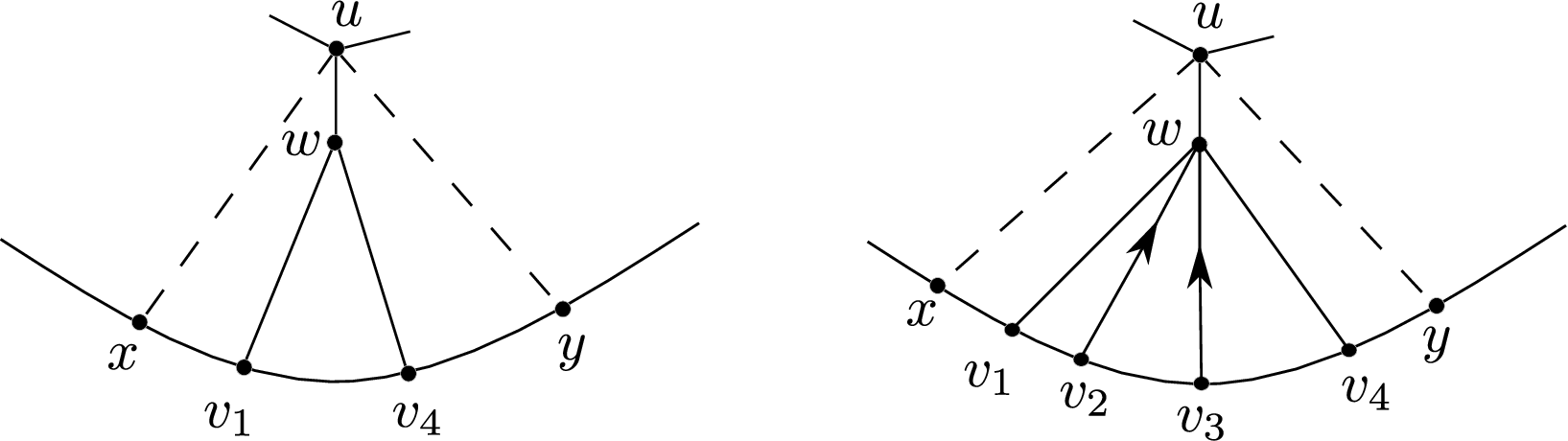}
\centerline{Fig.9 (for $m=4$) ~The graphs $G_4$ and $G$.}
\end{figure}

Since $m$ is even, the directed rim path $v_1v_2\cdots v_m$ has odd length $m-1$, so every Eulerian subdigraph of $D$ not using $\{a_2, \cdots, a_{m-1}\}$ has the same parity as the corresponding Eulerian subdigraph of $D_4$. Furthermore, since $m-2$ is even, pairing $a_{2k}$ with $a_{2k+1}$ for $k=1, \cdots, \frac{m-2}{2}$ changes the path length by $\pm 1$ and yields a parity-reversing bijection between $\bigcup_{i=2}^{m-1}\mathcal{E}_e(D, a_i)$ and $\bigcup_{i=2}^{m-1}\mathcal{E}_o(D, a_i)$. Hence ${\rm diff}(D)={\rm diff}(D_4)\neq 0$.

For each case, we obtain an $AT$-orientation $D$ of $G$ in which the outdegree of each vertex is at most 2, and thus $AT(G)\leq 3$. This completes the proof of Lemma 4.3.
\end{proof}

Combining Part 1 (for $G\in \mathcal{A}$) and Part 2 (Lemmas 4.1, 4.2, and 4.3 for $G\notin \mathcal{A}$), we conclude that $AT(G)=4$ if $G\in \mathcal{A}$ and $AT(G)=3$ otherwise. This completes the proof of Theorem 2. \hfill $\Box$

%%%%%%%%%%%
%\begin{figure}[htbp]
%\centering
%\includegraphics[height=4cm, width=0.7\textwidth]{2p3p5.eps}
%\centerline{Fig.2  ~The matching book embedding of $P_3\times P_5$.}%
%\end{figure}
%%%%%%%%%%%

%\begin{lem}
%If $H=T\cup C_{n}$ is not a wheel, and $C_{n}$ is an odd cycle, then $AT(H)=3$.
%\end{lem}

%$(u, v)$, $(v_{1}, v_{n})$ and $(v_{k}, v_{k+1})$ are not contained in any directed cycles. Hence none of these arcs is contained in any Eulerian subdigraph of $D$. Therefore
%$$|\mathcal{E}_{e}(D)|=|\mathcal{E}_{e}(D_{1})|\times |\mathcal{E}_{e}(D_{2})|+ |\mathcal{E}_{o}(D_{1})|\times |\mathcal{E}_{o}(D_{2})|,$$
%$$|\mathcal{E}_{o}(D)|=|\mathcal{E}_{e}(D_{1})|\times |\mathcal{E}_{o}(D_{2})|+ |\mathcal{E}_{o}(D_{1})|\times |\mathcal{E}_{e}(D_{2})|.$$
%diff$(D)$ = diff$(D_{1}) \times$ diff$(D_{2})\neq 0$, $D$ is an $AT$-orientation of $G$. Obvious that the followings hold:\\
%(1) $D$ is an $AT$-orientation.\\
%(2) $d^{+}_{D}(v)\leq2$ for each inner vertex $v$.\\
%(3) $d^{+}_{D}(v)\leq2$ for each outer vertex $v$.\\
%Hence $AT(H)\leq 3$.

%
%The result is established.

%According to $\chi(G)\leq ch(G)\leq \chi_{p}(G)\leq AT(G)$, we have the following result:
%$$\chi(H)=ch(H)=\chi_{P}(H)=AT(H)=\left\{\begin{array}{l}
%4,\ \mbox{if\ H\ is\ a\ wheel\ of\ even\ order};\\
%3,\ \mbox{otherwise}.\\
%\end{array} \right.$$

%%%%%%%%%%%
%\begin{figure}[htbp]
%\centering
%\includegraphics[height=4cm, width=0.7\textwidth]{3p3p5.eps}
%\centerline{Fig.3 ~The matching book embedding of $P_3\boxtimes P_5$.}
%\end{figure}
%%%%%%%%%%%

%%%%%%%%%%%%5

%%%%%%%%%%%%55


\begin{thebibliography}{99}

%\bibliographystyle{abbrv}

%\makeatletter

%\renewcommand\@biblabel[1]{#1}

%\makeatother


\bibitem{TR1995}%1
Jensen, T.R. and Toft, B., 1995. \emph{Graph coloring problems}. John Wiley $\&$ Sons, inc., New York. %pp.NY10158-0012.

\bibitem{JD2010}%3
Hladk\'{y}, J., Kr\'{a}l, D. and Schauz, U., 2010. Brooks' Theorem via the Alon-Tarsi Theorem. \emph{Discrete Mathematics}, \emph{310}(23), pp.3426-3428.

\bibitem{ESG2023}%4
Culver, E. and Hartke, S.G., 2023. Relation between the correspondence chromatic number and the Alon-Tarsi number. \emph{Discrete Mathematics},  \emph{346}(6), pp.113347.

\bibitem{XD2017}%7
Zhu, X.D., 2017. The Alon-Tarsi number of planar graph. \emph{Journal of Combinatorial Theory, Series B}, \emph{134}, pp.354$-$358.

\bibitem{JXD2020}%8
Grytczuk, J. and Zhu, X.D., 2020. The Alon-Tarsi number of planar graphs minus a matching. \emph{Journal of Combinatorial Theory, Series B}, \emph{145}, pp.511$-$520.

\bibitem{RS2019}%9
Kim, R., Kim, S. and Zhu, X.D., 2019. The Alon-Tarsi number of subgraphs of a planar graph. \emph{arXiv preprint arXiv:1906.01506}.

\bibitem{ZG2021}%10
Li, Z.G., Ye, Q. and Shao, Z.L., 2023. The Alon-Tarsi number of Halin graphs. \emph{Applied Mathematics-A Journal of Chinese Universities, A}, \emph{38}(03), pp.373-378.
%DOI:10.13299/j.cnki.amjcu.002278.

\bibitem{YQ1}%16
Li, Z.G., Ye, Q. and Shao, Z.L., 2023. The Alon-Tarsi number of two kinds of planar graphs. \emph{Electronic Journal of Graph Theory and Applications (EJGTA)}, \emph{11}(1), pp.155-164.

\bibitem{TAK2022}%6
Abe, T., Kim, S.J. and Ozeki, K., 2022. The Alon-Tarsi number of $K_5$-minor-free graphs. \emph{Discrete Mathematics}, \emph{345}(4), pp.112764.

\bibitem{ZG2023}%5
Li, Z.G., Shao, Z.L., Petrov, F. and Gordeev, A., 2023. The Alon-Tarsi number of toroidal grids. \emph{European Journal of Combinatorics}, \emph{111}, pp.103697.

\bibitem{YQ2}%17
Li, Z.G., Ye, Q. and Shao, Z.L., 2023. The Alon-Tarsi number of the Mycielski graphs. \emph{Journal of Applied Mathematics and Computation (JAMC)}, \emph{7}(1), pp.162-166.

\bibitem{XY2012}%12
Meng, X.Y., Guo, J.H. and Su, B.T., 2012. Incidence coloring of pseudo-Halin graphs. \emph{Discrete Mathematics}, \emph{312}(22), pp.3276-3282.

\bibitem{NM1992}%2
Alon, N. and Tarsi, M., 1992. Colorings and orientations of graphs. \emph{Combinatorica}, \emph{12}(2), pp.125$-$134.

\bibitem{LZ2001}%13
Liu, L.Z. and Zhang, Z.F., 2001. On the Structure Properties and Chromatics of pseudo-Halin Graphs. \emph{Journal of Lanzhou Railway Institute}, \emph{20}(04), pp.105-107.

\bibitem{NS2023}%14
Ghanbari, N. and Alikhani, S., 2023. Strong domination number of Haj\'{o}s sum and vertex-sum of two graphs. \emph{arXiv preprint arXiv:2302.01126}.

\bibitem{HX1994}%11
Li, H.X., Zhang, Z.F. and Zhang, J.X., 1994. On the colouring of Halin graphs. \emph{Journal of Shanghai Tiedao University (Natural Science Edition)}, \emph{15}(1), pp.19$-$24.


%\bibitem{D2011}%15
%Hefetz, D., 2011. On two generalizations of the Alon-Tarsi polynomial method. \emph{Journal of Combinatorial Theory, Series B}, 101(6), pp.403-414.

%\bibitem{XDR2021}%18
%Zhu, X.D. and Balakrishnan, R., 2022. \emph{Combinatorial Nullstellensatz: With Applications to Graph Colouring}. CRC Press, Boca Raton.% pp.FL33487-242.


%\bibitem{EN}%7:
%J. A. Bondy, U. S. R. Murty. Graph Theory[M]. London: Springer, 2008.





\end{thebibliography}
\end{document}